# On the Mechanism of Roe-type Schemes for All-Speed Flows

Xue-song Li [*], Chun-wei Gu

*Key Laboratory for Thermal Science and Power Engineering of Ministry of Education, Department of Thermal Engineering, Tsinghua University, Beijing 100084, PR China*

**Abstract:** In recent years, Roe-type schemes based on different ideas have been developed for all-speed flows, such as the preconditioned Roe, the All-Speed Roe, Thornber's modified Roe and the LM-Roe schemes. This work explores why these schemes succeed or fail with the accuracy and checkerboard problems. Comparison and analysis show that the accuracy and checkerboard problems are caused by the order of the sound speed being too large and too small in the coefficients of the velocity-derivative and pressure-derivative dissipation terms, respectively. These problems can be resolved by choosing coefficients with zero-order sound speed. In addition, to avoid the negative effects of the global cut-off strategy on accuracy while maintaining computational stability, the sound speed terms in the numerator of the coefficients can be determined by local variables, while those in the denominator remain the global cut-off. Two novel schemes are proposed as examples to demonstrate how these ideas can be applied to construct more satisfactory schemes for all-speed flows. Asymptotic analysis and numerical experiments support the theoretical analysis and the rules obtained in the work.

## 1. Introduction

[*] Corresponding author. Tel.: 0086-10-62794617; fax: 0086-10-62795946
E-mail address: xs-li@mail.tsinghua.edu.cn (X.-S. Li).

The Roe scheme is a classical and popular shock-capturing scheme that is widely used in the computation of compressible flows. For nearly-incompressible flows, however, it produces unphysical discrete results because of the incorrect scaling of pressure fluctuations, as shown by asymptotic analysis [1]. To overcome this limitation, preconditioning methods have been developed and applied to the Roe scheme in recent decades, resulting in a new scheme named the preconditioned Roe scheme [2,3], which shows the correct asymptotic behavior at the low Mach number limit [1].

Like other preconditioned schemes, however, the preconditioned Roe scheme suffers from the global cut-off problem [3]. Without adopting the global cut-off strategy, preconditioned schemes are unstable unless the time step is extremely small, especially for viscous flows. Adopting it, however, limits the ability to accurately simulate low- and high-Mach number mixed flows. For example, under extreme conditions, in which shock waves and incompressible flows coexist in the same flow region, the global cut-off will reduce the preconditioned Roe scheme to the Roe scheme in the whole flow region; calculation of the incompressible region cannot then benefit from the preconditioning and will suffer from the accuracy problem. From this point of view, the preconditioned Roe scheme cannot be regarded as a universal scheme for all-speed flows.

In recent years, improved Roe-type schemes have been developed that only need local parameters, without the global cut-off strategy for all-speed flows. For example, the All-Speed Roe scheme [4-7] proposes only modifying the non-linear eigenvalues; Thornber [8] finds the modification of the right eigenvector matrix; and the fix recommended in [9-11] is based on the idea [12] that the accuracy problem is attributed to the velocity jump normal to the cell interface. Although based on different ideas, all these methods achieve the same purpose. This fact motivates us to understand the common underlying mechanism of these schemes and then find the general rules needed to construct satisfactory schemes.

The outline of this paper is as follows. Section 2 presents the governing equations, the general form of the schemes and a discussion of the central terms. Section 3 briefly reviews the Roe, the preconditioned Roe, the All-Speed Roe, Thornber's modified Roe (T-Roe), and the low-Mach number fix for Roe (LM-Roe) schemes. Section 4 discusses the mechanism underlying these schemes and conducts an asymptotic analysis to demonstrate it. Section 5 proposes two novel Roe-type schemes that show how the mechanism can be used to construct more satisfactory schemes. Section 6 presents the numerical experiments to support the theoretical analysis and prediction. Finally, section 7 closes the paper with some concluding remarks and a summary of the general rules for improving schemes.

## 2. Governing Equations and the Central Term of Schemes

### 2.1 Governing Equations

For simplicity, the two-dimensional Euler compressible equations are written as

$$\frac{\partial \boldsymbol{Q}}{\partial t} + \frac{\partial \boldsymbol{F}}{\partial x} + \frac{\partial \boldsymbol{G}}{\partial y} = 0, \tag{1}$$

where: $Q = \begin{bmatrix} \rho \\ \rho u \\ \rho v \\ \rho E \end{bmatrix}$ is the vector of conversation variables, $F = \begin{bmatrix} \rho u \\ \rho u^2 + p \\ \rho uv \\ u(\rho E + p) \end{bmatrix}$ and

$G = \begin{bmatrix} \rho v \\ \rho uv \\ \rho v^2 + p \\ v(\rho E + p) \end{bmatrix}$ are the vectors of Euler fluxes, $\rho$ is the fluid density, $p$ is the pressure, $E$ is the total energy, and $u, v$ are the velocity components in the Cartesian coordinates $(x, y)$, respectively.

To avoid the checkerboard pressure-velocity decoupling problem in the computation of incompressible flows, the concept of interface velocity, i.e. the velocity on the interface between the adjacent control volumes, can be introduced using the time-marching momentum interpolation method (MIM) [7], and then the governing equations of Eq. (1) become:

$$\frac{\partial Q}{\partial t} + \frac{\partial F^{MIM}}{\partial x} + \frac{\partial G^{MIM}}{\partial y} = 0, \qquad (2)$$

where $F^{MIM} = u_f \begin{bmatrix} \rho \\ \rho u \\ \rho v \\ \rho H \end{bmatrix} + \begin{bmatrix} 0 \\ p \\ 0 \\ 0 \end{bmatrix}$, $G^{MIM} = v_f \begin{bmatrix} \rho \\ \rho u \\ \rho v \\ \rho H \end{bmatrix} + \begin{bmatrix} 0 \\ 0 \\ p \\ 0 \end{bmatrix}$, $H = E + \frac{p}{\rho}$, and $u_f$ and $v_f$

are the components of the interface velocity.

## 2.2 General Form of Schemes and Discussion of the Central Term

Many schemes that include the Roe scheme can be generalized as the sum of a central term $\tilde{F}_c$ and a numerical dissipation term $\tilde{F}_d$:

$$\tilde{F} = \tilde{F}_c + \tilde{F}_d. \qquad (3)$$

For simplicity, the equations of the schemes are given only for the *i*-direction with the index *j* omitted.

In general, the central term $\tilde{F}_c$ is obtained by averaging the fluxes as:

$$\tilde{F}_{c,i+\frac{1}{2}} = \frac{1}{2}\left(\bar{F}_i + \bar{F}_{i+1}\right) = \frac{1}{2}\left\{ \begin{bmatrix} \rho U \\ \rho U u + n_x p \\ \rho U v + n_y p \\ \rho U H \end{bmatrix}_i + \begin{bmatrix} \rho U \\ \rho U u + n_x p \\ \rho U v + n_y p \\ \rho U H \end{bmatrix}_{i+1} \right\}, \tag{4}$$

where $U = n_x u + n_y v$ is the velocity normal to the cell interface, and $n_x$ and $n_y$ are the components of the face-normal vector, respectively.

As discussed in Ref. [5-7], some preconditioned schemes with Eq. (4) suffer from the checkerboard problem. To avoid this problem, the governing equations of Eq. (2) can be applied to the central term as:

$$\tilde{F}_{c,i+\frac{1}{2}} = \frac{1}{2}U_f \left\{ \begin{bmatrix} \rho \\ \rho u \\ \rho v \\ \rho H \end{bmatrix}_i + \begin{bmatrix} \rho \\ \rho u \\ \rho v \\ \rho H \end{bmatrix}_{i+1} \right\} + \frac{1}{2}\left\{ \begin{bmatrix} 0 \\ n_x p \\ n_y p \\ 0 \end{bmatrix}_i + \begin{bmatrix} 0 \\ n_x p \\ n_y p \\ 0 \end{bmatrix}_{i+1} \right\}, \tag{5}$$

where $U_f = n_x u_f + n_y v_f$ is the normal interface velocity, which can be obtained as follows:

$$\left(U_f\right)^n_{i+\frac{1}{2}} = 0.5\left(U_i^n + U_{i+1}^n\right) + \delta U_f, \tag{6}$$

where $\delta U_f$ is the modification of $U_f$. When

$$\delta U_f = 0, \tag{7}$$

Eq. (6) becomes the second-order central interpolation and the corresponding fluxes defined in Eq. (5) will be very close to those in Eq. (4) at the low Mach number limit. To ensure the velocity-pressure coupling, Ref. [13, 14] propose adding a first pressure-derivative smoothing term to the interface velocity, and this idea is adopted by the All-Speed Roe scheme as [4-6]:

$$\delta U_f = -\frac{c_2}{\rho^* u^*}\left(p_{i+1}^n - p_i^n\right), \tag{8}$$

where $\rho^*$ and $u^*$ are the reference density and velocity respectively, and $c_2$ is a constant. $c_2$ should be as small as possible for accuracy only if the checkerboard decoupling is suppressed, and too large a value of $c_2$ can lead to computational divergence, as analyzed in section 4.4. Although the optimum value of $c_2$ depends on specific problems, a value of 0.04 is recommended; this shows good applicability, as demonstrated in Chapter 6.

Ref. [7] proposes a time-marching momentum interpolation method (MIM) for steady and unsteady flows. The following is the equation for the method adapted for steady simulations with a large time-step:

$$\delta U_f = n_x \Delta t \left[\frac{\left(\frac{\partial p}{\partial x}\right)_i}{2\rho_i} + \frac{\left(\frac{\partial p}{\partial x}\right)_{i+1}}{2\rho_{i+1}} - \frac{\left(\frac{\partial p}{\partial x}\right)_{i+1/2}}{\rho_{i+1/2}}\right]^{n-1} + n_y \Delta t \left[\frac{\left(\frac{\partial p}{\partial y}\right)_i}{2\rho_i} + \frac{\left(\frac{\partial p}{\partial y}\right)_{i+1}}{2\rho_{i+1}} - \frac{\left(\frac{\partial p}{\partial y}\right)_{i+1/2}}{\rho_{i+1/2}}\right]^{n-1},$$

$$\tag{9}$$

which introduces a third pressure-derivative smoothing term into the interface velocity, leading to higher accuracy than Eq. (8) [7]. To calculate $\left(\frac{\partial p}{\partial x}\right)_{i+1/2}$, the algorithm of the Gauss integral needs information on all neighbour nodes but does not make codes more difficult, especially for viscous calculation, because the same algorithm is needed to calculate the velocity gradient at $x_{i+1/2}$.

Although the above equations are regarded as the central term, in fact Eqs. (8) and (9) introduce the pressure-difference type numerical dissipation and are similar to $\delta U$ in Eq. (30), the dissipation term discussed in Chapter 4. Therefore, $\delta U_f$ is actually part of $\delta U$.

## 3. The Numerical Dissipation Term of the Schemes

### 3.1 The Roe Scheme

For the "classical" Roe scheme, the first-order numerical dissipation term can be expressed as:

$$\tilde{\boldsymbol{F}}^{Roe}_{d,i+\frac{1}{2}} = -\frac{1}{2} \boldsymbol{R}^{Roe}_{i+\frac{1}{2}} \left| \boldsymbol{\Lambda}^{Roe}_{i+\frac{1}{2}} \right| \left( \boldsymbol{R}^{Roe}_{i+\frac{1}{2}} \right)^{-1} (\boldsymbol{Q}_{i+1} - \boldsymbol{Q}_i), \tag{10}$$

where $\boldsymbol{R}^{Roe}$ is the right eigenvector matrix of $\dfrac{\partial \boldsymbol{F}}{\partial \boldsymbol{Q}}$:

$$\boldsymbol{R}^{Roe} = \begin{bmatrix} 1 & 0 & 1 & 1 \\ u & -n_y & u-n_x c & u+n_x c \\ v & n_x & v-n_y c & v+n_y c \\ \frac{1}{2}(u^2+v^2) & n_x v - n_y u & H-cU & H+cU \end{bmatrix}, \tag{11}$$

where $U = n_x u + n_y v$ is the velocity normal to the cell interface, $c$ is the sound speed, and $n_x$ and $n_y$ are the components of the face normal vector.

In addition, $\boldsymbol{\Lambda}^{Roe}$ is the diagonal matrix consisting of the relevant eigenvalues:

$$\boldsymbol{\Lambda}^{Roe} = \begin{bmatrix} U & & & \\ & U & & \\ & & U-c & \\ & & & U+c \end{bmatrix}. \tag{12}$$

### 3.2 The Preconditioned Roe Scheme

The preconditioned Roe scheme (P-Roe) is expressed as:

$$\tilde{F}_{d,i+\frac{1}{2}}^{P-Roe} = -\frac{1}{2}\Gamma_{i+\frac{1}{2}}^{-1} R_{i+\frac{1}{2}}^{P-Roe} \left|\Lambda_{i+\frac{1}{2}}^{P-Roe}\right| \left(R_{i+\frac{1}{2}}^{P-Roe}\right)^{-1} (Q_{i+1} - Q_i), \tag{13}$$

where $\Gamma$ is the preconditioner based on the conservative variables, as proposed in Ref. [1-3]:

$$\Gamma = \frac{\partial Q}{\partial W} \begin{bmatrix} \theta & & & \\ & 1 & & \\ & & 1 & \\ & & & 1 \end{bmatrix} \frac{\partial W}{\partial Q}, \tag{14}$$

where $W$ is a vector of the primitive variables $[p\ u\ v\ S]^T$. When $\theta = 1$, P-Roe becomes the "classical" Roe scheme.

In addition, $R^{P-Roe}$ is the right eigenvector matrix of $\Gamma \frac{\partial F}{\partial Q}$ and $\Lambda^{P-Roe}$ is the diagonal matrix including the relevant eigenvalues.

$$R^{P-Roe} = \begin{bmatrix} 1 & 0 & 1 & 1 \\ u & -n_y & u - \frac{\theta-1}{2\theta}n_x U - n_x \frac{\tilde{c}}{\theta} & u - \frac{\theta-1}{2\theta}n_x U + n_x \frac{\tilde{c}}{\theta} \\ v & n_x & v - \frac{\theta-1}{2\theta}n_y U - n_y \frac{\tilde{c}}{\theta} & v - \frac{\theta-1}{2\theta}n_y U + n_y \frac{\tilde{c}}{\theta} \\ \frac{1}{2}(u^2+v^2) & n_x v - n_y u & H - \frac{\theta-1}{2\theta}U^2 - U\frac{\tilde{c}}{\theta} & H - \frac{\theta-1}{2\theta}U^2 + U\frac{\tilde{c}}{\theta} \end{bmatrix}, \tag{15}$$

and

$$\Lambda^{P-Roe} = \begin{bmatrix} U & & & \\ & U & & \\ & & \tilde{U}-\tilde{c} & \\ & & & \tilde{U}+\tilde{c} \end{bmatrix}, \tag{16}$$

where $\tilde{U} = \frac{1}{2}(1+\theta)U$, $\tilde{c} = \frac{1}{2}\sqrt{4c^2\theta + (1-\theta)^2 U^2}$, and $\theta$ is the key factor for the preconditioning technique. In theory, the value of the key factor $\theta$ should be related to the local Mach number. In practice, however, $\theta$ is cut off by the global Mach number to avoid computational instability:

$$\theta = \min\left[\max\left(KM_{ref}^2, M^2\right), 1\right], \tag{17}$$

where the constant $K$ typically equals 1, and the reference Mach number $M_{ref}$ is a global parameter, which may be the inlet Mach number, or the average or maximum Mach number over the flow. This means that the accuracy will be high enough in high-speed flow regions, such as the main flows, but deteriorates in low-speed flow regions, such as the boundary layers.

For nearly-incompressible flows, the pseudo-acoustic speed with the global cut-off is $\tilde{c} = \max\left(\frac{\sqrt{5}}{2}M, M_{ref}\right)c$. For clear comparison and unified expression, it can be redefined as:

$$\tilde{c} = \tilde{f}\left(M, M_{ref}\right)c, \tag{18}$$

where

$$\tilde{f}\left(M, M_{ref}\right) = \min\left[\max\left(M, M_{ref}\right), 1\right]. \tag{19}$$

The pseudo-acoustic speed without the cut-off $c'$ can be defined according to the local Mach number as:

$$c' = f(M)c, \tag{20}$$

where

$$f(M) = \min(M, 1). \tag{21}$$

### 3.3 The All-Speed Roe Scheme

In preconditioning techniques, a global cut-off strategy is adopted to avoid computational instability, which is attributed in Ref. [4] to the structure of $\frac{1}{\theta}$ or $\frac{1}{\tilde{c}}$, because it produces the term $\frac{1}{U}$, which can greatly magnify the velocity fluctuation when the velocity is small but with a steep gradient. To avoid the global cut-off, therefore, Ref. [4] proposes two procedures: (a) the governing equations are regarded as the aggregated scalar equations, although they are usually treated as a coupling system under the compressible condition; and (b) the discretization of the pressure-derivative term is no longer consistent with that of the convective terms and independently adopts a centred scheme. Then, extending the Roe scheme to the governing system defined by these procedures leads to the Low-Speed Roe scheme that is only applicable to low-Mach-number flows:

$$\tilde{F}_{d,i+\frac{1}{2}}^{L-Roe} = -\frac{1}{2}\left|U_{i+\frac{1}{2}}\right|(Q_{i+1} - Q_i) \tag{22}$$

An idea [5, 6] is further proposed to unify the Low-Speed Roe scheme and shock-capturing scheme with a function of the local Mach number. Based on this idea, the All-Speed Roe scheme is derived as:

$$\tilde{F}_{d,i+\frac{1}{2}}^{A-Roe} = -\frac{1}{2}R_{i+\frac{1}{2}}^{Roe}\left|\Lambda_{i+\frac{1}{2}}^{A-Roe}\right|\left(R_{i+\frac{1}{2}}^{Roe}\right)^{-1}(Q_{i+1} - Q_i), \tag{23}$$

where $R^{Roe}$ has the same right eigenvector matrix of $\frac{\partial F}{\partial Q}$ as for the "classical" Roe scheme (11). However, the elements of the diagonal matrix $\Lambda^{A-Roe}$ are different from those in Eq. (12):

$$\Lambda^{A-Roe} = \begin{bmatrix} U & & & \\ & U & & \\ & & U-c' & \\ & & & U+c' \end{bmatrix}, \tag{24}$$

where $c' = f(M)c$.

It can be seen that the All-Speed Roe scheme simply multiplies the acoustic term $c$ in the eigenvalues of the Roe scheme by the factor $f(M)$ to remedy the accuracy problem. The rationality for this procedure has been verified theoretically with an asymptotic analysis [6] and Hodge decomposition [15].

The function $f(M)$ is only related to the local Mach number and its choice is not unique [5, 6]. In the present work, Eq. (21) can be adopted for comparison.

### 3.4 The T-Roe Scheme

Thornber et al. [8] find that the following modification of the right eigenvector matrix can also lead to an All-Speed Roe-type scheme without the limitation of the global cut-off:

$$\tilde{F}_{d,i+\frac{1}{2}}^{T-Roe} = -\frac{1}{2} R_{i+\frac{1}{2}}^{T-Roe} \left| \Lambda_{i+\frac{1}{2}}^{Roe} \right| \left( R_{i+\frac{1}{2}}^{Roe} \right)^{-1} (Q_{i+1} - Q_i), \tag{25}$$

where $R^{T-Roe}$ is the modified $R^{Roe}$, in which the acoustic speed terms have been multiplied by the function $f(M)$:

$$R^{T-Roe} = \begin{bmatrix} 1 & 0 & 1 & 1 \\ u & -n_y & u - n_x c' & u + n_x c' \\ v & n_x & v - n_y c' & v + n_y c' \\ \frac{1}{2}(u^2 + v^2) & n_x v - n_y u & H - Uc' & H + Uc' \end{bmatrix}, \tag{26}$$

where $c'$ can adopt the form defined in Eqs. (20) and (21).

It can be seen that in the original version [8], the components $R^{T-Roe}(4,3)$ and $R^{T-Roe}(4,4)$ are equal to those in $R^{Roe}$ without the modification, indicating that the energy equation is not necessarily modified. This is not surprising because the term $cU$ can be disregarded in comparison to the term $H$ of the order $\mathrm{O}(c^2)$ for the low Mach number limit, which makes the effort of decreasing this term trivial. The aim of modifying the energy equation in Eq. (26) is only for the convenience of uniformity of expression and comparison.

### 3.5 The LM-Roe Scheme

In the Roe scheme, Eq. (10) $\left(\boldsymbol{R}^{Roe}_{i+\frac{1}{2}}\right)^{-1}(\boldsymbol{Q}_{i+1}-\boldsymbol{Q}_i)$ can be rewritten as follows, without considering the conservative characteristics:

$$\boldsymbol{\alpha}^{Roe}_{i+\frac{1}{2}} = \left(\boldsymbol{R}^{Roe}_{i+\frac{1}{2}}\right)^{-1}(\boldsymbol{Q}_{i+1}-\boldsymbol{Q}_i) = \begin{bmatrix} \Delta\rho - \dfrac{\Delta p}{c^2} \\ \rho\Delta V \\ \dfrac{\Delta p}{2c^2} - \dfrac{\rho}{2c}\Delta U \\ \dfrac{\Delta p}{2c^2} + \dfrac{\rho}{2c}\Delta U \end{bmatrix}, \qquad (27)$$

where $V = n_x v - n_y u$.

As mentioned above, the accuracy problem of the Roe scheme is attributed to $\Delta U$ [12], the velocity jump normal to the cell interface. Based on this idea, an improved reconstruction method is proposed [11]. Fillion [9] and Rieper [10] also propose a low Mach number fix, in which $\Delta U$ is multiplied by the function $f(M)$ and Eq. (27) is rewritten as:

$$\boldsymbol{\alpha}^{LM-Roe}_{i+\frac{1}{2}} = \begin{bmatrix} \Delta\rho - \dfrac{\Delta p}{c^2} \\ \rho\Delta V \\ \dfrac{\Delta p}{2c^2} - \dfrac{\rho}{2c}f(M)\Delta U \\ \dfrac{\Delta p}{2c^2} + \dfrac{\rho}{2c}f(M)\Delta U \end{bmatrix}, \qquad (28)$$

where $f(M)$ can also adopt the form of Eq. (21).

Therefore, the LM-Roe scheme can be expressed as:

$$\tilde{\boldsymbol{F}}_{d,i+\frac{1}{2}}^{LM-Roe} = -\frac{1}{2}\boldsymbol{R}_{i+\frac{1}{2}}^{Roe}\left|\boldsymbol{\Lambda}_{i+\frac{1}{2}}^{Roe}\right|\boldsymbol{\alpha}_{i+\frac{1}{2}}^{LM-Roe} \tag{29}$$

## 4. Analysis of the Mechanism of the Roe-type Schemes for All-Speed Flows

### 4.1 A General Form for the Mechanism Analysis

All of these Roe-type schemes are successful for all-speed flows regardless of which ideas they are based on, with or without the global cut-off. To understand the underlying common mechanism, the numerical dissipation terms in Chapter 3 are expanded and rewritten as a unified form. Inspired by Ref. [2], a form that is valid for all these schemes with the low Mach number limit is proposed as:

$$\tilde{\boldsymbol{F}}_d = -\frac{1}{2}\left\{|U|\begin{bmatrix}\Delta\rho \\ \Delta(\rho u) \\ \Delta(\rho v) \\ \Delta(\rho E)\end{bmatrix} + \delta U\begin{bmatrix}\rho \\ \rho u \\ \rho v \\ \rho H\end{bmatrix} + \delta p\begin{bmatrix}0 \\ n_x \\ n_y \\ U\end{bmatrix}\right\}, \tag{30}$$

which includes three terms on the right side. The first term is the basic upwind dissipation, which is consistent with Eq. (22), the numerical dissipation term of the Low-Speed Roe scheme. The second term has a role similar to the MIM fluxes defined in Eqs. (8) and (9). The third term is a modification to the interface pressure, which will decide the accuracy, as described in detail in section 4.3.

### 4.2 The Mechanism Analysis of the Modification to the Normal Interface Velocity

Specific expressions of $\delta U$ for the various schemes are given as follows.

For the Roe scheme,

$$\delta U = (c - |U|)\frac{\Delta p}{\rho c^2} + \frac{U}{c}\Delta U. \tag{31}$$

For the preconditioned Roe scheme,

$$\delta U = \left[\tilde{c} - \frac{1-\theta}{2}U\frac{\tilde{U}}{\tilde{c}} - \theta|U|\right]\frac{\Delta p}{\rho\theta c^2} + \frac{\tilde{U}}{\tilde{c}}\Delta U. \tag{32}$$

For the All-Speed Roe scheme,

$$\delta U = (c' - |U|)\frac{\Delta p}{\rho c^2} + \frac{U}{c}\Delta U. \tag{33}$$

For the T-Roe scheme,

$$\delta U = (c - |U|)\frac{\Delta p}{\rho c^2} + \frac{U}{c}\Delta U. \tag{34}$$

For the LM-Roe scheme

$$\delta U = (c - |U|)\frac{\Delta p}{\rho c^2} + \frac{U}{c}f(M)\Delta U. \tag{35}$$

In all these equations, $\tilde{c} = \tilde{f}(M, M_{ref})c$, and $c' = f(M)c$, as defined in Eq. (18) and Eq. (20), which ensures $\tilde{c} \geq |\tilde{U}|$ and $c' \geq |U|$.

As discussed in section 4.1, the pressure-derivative term in $\delta U$ is important in suppressing the checkerboard problem because it has a similar effect to the time-marching MIM, and the velocity-derivative term $\frac{U}{c}\Delta U$ can be disregarded both from the point of view of the mechanism and the amount of dissipation.

By comparing Eqs. (31) – (35), it can be seen that the equations for the Roe scheme (31), the T-Roe scheme (34), and the LM-Roe scheme (35) are equivalent if the velocity-derivative term is not taken into account. This suggests that $\delta U$ in the Roe scheme does not need modification, at least from the viewpoint of accuracy.

Although $\delta U$ has little effect on the accuracy of each scheme, the pressure-derivative term has different effects on the checkerboard problem for different schemes. The coefficient of the pressure-derivative term reflects the ability to suppress the checkerboard problem. The orders of this coefficient are $\mathrm{O}(c^{-1})$ for the T-Roe and LM-Roe schemes Eqs. (34) and (35), $\mathrm{O}(c^0)$ for the preconditioned Roe scheme Eq. (32), and $\mathrm{O}(c^{-2})$ for the All-Speed Roe scheme Eq. (33). Thus, the preconditioned Roe scheme is the best for avoiding checkerboard, the All-Speed Roe scheme is the worst, and the T-Roe and the LM-Roe schemes are in-between. This estimation accords with the analyses given in Refs. [1,6,7,10] and will be discussed further in section 4.4 with an asymptotic analysis. In summary, a pressure-derivative term of the order $\mathrm{O}(c^0)$ in $\delta U$ is required and can be provided by Eq. (8) or Eq. (9).

### 4.3 Mechanism Analysis of the Modification to the Interface Pressure

Specific expressions of $\delta p$ for the various schemes are given as follows.

For the Roe scheme,

$$\delta p = \frac{U}{c}\Delta p + (c - |U|)\rho\Delta U . \tag{36}$$

For the preconditioned Roe scheme,

$$\delta p = \frac{\tilde{U}}{\tilde{c}}\Delta p + \left(\tilde{c} - |U| + \frac{1-\theta}{2}U\frac{\tilde{U}}{\tilde{c}}\right)\rho\Delta U . \tag{37}$$

For the All-Speed Roe scheme,

$$\delta p = \frac{U}{c}\Delta p + (c' - |U|)\rho\Delta U . \tag{38}$$

For the T-Roe scheme,

$$\delta p = f(M)\frac{U}{c}\Delta p + (c' - |U|)\rho\Delta U . \tag{39}$$

For the LM-Roe scheme,

$$\delta p = \frac{U}{c}\Delta p + (c' - |U|)\rho\Delta U . \tag{40}$$

A common feature of Eqs. (37)-(40) is that the coefficient $c$ of the term $c\rho\Delta U$ in Eq. (36) is replaced by another coefficient of the order $\mathrm{O}(c^0)$, such as $\tilde{c}$ or $c'$.

The effect of $\delta p$ on the accuracy is limited to the term $c\rho\Delta U$ in the momentum equation for the Roe scheme, because it disappears in the continuity equation and can be disregarded in the energy equation in comparison to $H$ of the order $\mathrm{O}(c^2)$, as discussed in section 3.4. The problem, therefore, can be solved by modifying only this term.

It should also be noted that Eqs. (38)-(40) all adopt $c'$, decided by the local variables, to modify the term $c\rho\Delta U$, whereas Eq. (37) adopts $\tilde{c}$ with the global cut-off. This further validates the analysis given in section 3.3, which attributes the global cut-off to the term $\frac{1}{\tilde{c}}$ or $\frac{1}{U}$.

**4.4 Asymptotic Analysis**

As discussed above, the accuracy problem in the Roe scheme is caused by the coefficient of the term $\Delta U$ in $\delta p$ for the momentum equation, and the checkerboard problem is caused by the coefficient of the term $\Delta p$ in $\delta U$. To further verify these conclusions, an asymptotic analysis, which is widely utilized in [1, 6, 7, 10] to investigate the behavior of continuous or discrete flows with a low Mach number, is performed as follows.

Eq. (30) can be recognized as the first basic upwind dissipation term plus a velocity-derivative term and a pressure-derivative term. Thus, the continuity equation becomes

$$\tilde{F}_{d,\rho} = -\frac{1}{2}\left[|U|\Delta\rho + g_\rho \Delta U + h_\rho \Delta p\right], \tag{41}$$

and the momentum equation

$$\tilde{F}_{d,\rho u} = -\frac{1}{2}\left[|U|\Delta(\rho \boldsymbol{u}) + g_{\rho u}\Delta \boldsymbol{u} + h_{\rho u}\Delta p\right], \tag{42}$$

where the coefficients $g_\rho$, $h_\rho$, and $h_{\rho u}$ are of the order $O(c^{-1})$ and $g_{\rho u}$ is of the order $O(c^1)$ for the Roe scheme.

To analyze a steady low-Mach-number flow, all non-dimensional variables are asymptotically expanded into powers of the reference Mach number $M_*$:

$$\tilde{\phi} = \tilde{\phi}^0 + M_*\tilde{\phi}^1 + M_*^2\tilde{\phi}^2 + M_*^3\tilde{\phi}^3 + \cdots,$$

where $\phi$ represents one of the fluid variables, $\rho$, $u$, $v$, $E$, or $p$.

If $g_{\rho u}$ is modified, from the order of $O(c^1)$ in the Roe scheme to the order of $O(c^0)$, then from the asymptotically expanded momentum equation (42) the following terms with an equal power of $M_*$ can be collected:

$$p^0_{i-1,j} - p^0_{i+1,j} = 0, \tag{43}$$

$$p^0_{i,j-1} - p^0_{i,j+1} = 0, \tag{44}$$

$$p^1_{i-1,j} - p^1_{i+1,j} = 0, \tag{45}$$

$$p^1_{i,j-1} - p^1_{i,j+1} = 0. \tag{46}$$

Leaving aside the possibility of the checkerboard problem, Eqs. (43) – (46) result in $p^0_i = cte \ \forall \boldsymbol{i}$ and $p^1_i = cte \ \forall \boldsymbol{i}$, which satisfies the physical behavior of low-Mach-number flows. Whereas, if $g_{\rho u}$ retains the order of $O(c^1)$, Eqs. (45) and

(46) no longer hold, leading to the accuracy problem of the Roe scheme.

When $h_\rho$ is of the order $O(c^0)$, the following terms with an equal power of $M_*$ can be collected from the asymptotically expanded continuity equation (41):

$$\sum_{l \in v(i)} h_\rho^0 \Delta_{il} p^0 = 0, \tag{47}$$

$$\sum_{l \in v(i)} h_\rho^0 \Delta_{il} p^1 = 0, \tag{48}$$

which can suppress the possible checkerboard of $p^0$ and $p^1$, respectively. The numerical experiments in Chapter 6 show that the checkerboard of $p^2$ can also be avoided.

If $h_\rho$ is of the order $O(c^{-1})$, Eq. (48) no longer holds, and may lead to the checkerboard of $p^1$. As shown in Ref. [10], however, this checkerboard mode can be damped away by the numerical viscosity due to the maintained upwinding in the momentum equations. Unfortunately, the checkerboard of $p^2$ can also occur, as shown in Chapter 6. Furthermore, it is not difficult to deduce that $p^0$ suffers from a severe checkerboard problem if $h_\rho$ is of the order $O(c^{-2})$, as in the All-Speed Roe scheme with the central term Eq. (4). It is clear that $h_\rho$ with a higher order of $c$ can provide a stronger restriction on the checkerboard.

In addition, it should be noticed that although $h_\rho$ of the order $O(c^0)$ is necessary for suppressing the checkerboard problem, it will result in the following physical feature [6] being unsatisfied:

$$u_{i+1,j}^0 - u_{i-1,j}^0 + v_{i,j+1}^0 - v_{i,j-1}^0 = 0. \tag{49}$$

For example, in the preconditioned Roe scheme [6]:

$$u^0_{i+1,j} - u^0_{i-1,j} + v^0_{i,j+1} - v^0_{i,j-1} = -\frac{1}{\rho_l^0 E_l^0 + p_l^0} \sum_{l \in v(i)} \frac{h_{il}^0}{2\sqrt{Y_{il}^0}} U_{il}^0 \rho_{il}^0 \Delta_{il} U^0 + \frac{h_{il}^0}{\sqrt{Y_{il}^0}} \Delta_{il} p^2 ,$$

and in the All-Speed Roe scheme with Eq. (8) [7]:

$$u^0_{i+1,j} - u^0_{i-1,j} + v^0_{i,j+1} - v^0_{i,j-1} = c_2 (\Delta x)^2 \nabla^2 p^2 .$$

The latter also explains why attempting to increase the order of $h_\rho$ to $\mathrm{O}(c^1)$ will fail and why too large a value of $c_2$ in Eq. (8) will lead to divergence of computation, because there is too much deviation from the condition in Eq. (49).

A possible way of avoiding this difficulty is the time-marching MIM Eq. (9) described in section 2.2, which produces a spatially third-order error, and the zero velocity divergence condition (49) is satisfied from the numerical viewpoint if popular second-order schemes are used [7]:

$$u^0_{i+1,j} - u^0_{i-1,j} + v^0_{i,j+1} - v^0_{i,j-1} = \frac{\mathrm{CFL}}{\rho^0 u^0} (\Delta x)^4 \nabla^4 p^2 .$$

**5. Two Novel Roe-type Schemes Based on the Mechanisms of All-Speed Flows**

According to the analysis in Chapter 4, two rules can be summarized as follows.

1. The accuracy problem is controlled by the coefficients of the $\Delta U$ term, and the checkerboard problem by the coefficients of the $\Delta p$ term. Both are determined by the order of $c$ in the respective coefficients.

2. The problem of the global cut-off is due to the pseudo-sound speed in the denominator of the coefficients.

These two rules reveal the deeper mechanisms underlying the schemes for All-Speed flows, and suggest directions for constructing more satisfactory schemes.

Two novel Roe-type schemes are proposed here to validate the discoveries.

The first version is based on the first rule, which involves the least modification to the Roe scheme by replacing $c$ with $c'$ in the $\Delta U$ term of $\delta p$:

$$\delta U = \left(c - |U|\right)\frac{\Delta p}{\rho c^2} + \frac{U}{c}\Delta U, \tag{50}$$

$$\delta p = \frac{U}{c}\Delta p + \left(c' - |U|\right)\rho \Delta U. \tag{51}$$

Compared with the LM Roe scheme, it merely replaces $\frac{U}{c}f(M)\Delta U$ with $\frac{U}{c}\Delta U$ in $\delta U$. From a numerical viewpoint, there is little difference between this new scheme and the LM Roe scheme, because the term $\frac{U}{c}\Delta U$ is trivial at the low Mach number limit. However, they reflect different perspectives on the mechanisms of All-Speed schemes, with the accuracy problem attributed to the coefficient of the $\Delta U$ term and the $\Delta U$ term itself [9-12], respectively.

In more detail, although the opinion on the accuracy problem in the first rule looks similar to the opinion in Ref. [9-12], they are quite different in nature. The opinion in Ref. [9-12] indicates that all $\Delta U$ terms should be corrected by $f(M)$, because the accuracy problem is attributed to the jump in the normal velocity component. However, this opinion suffers from two problems. One is an inherent contradiction of the opinion that $\Delta U$ does not need to be corrected in the $|U|\rho\Delta U$ term in $\delta p$ of the LM-Roe scheme. The other is the fact that other schemes discussed in the paper are successful even though they do not correct all $\Delta U$ terms. Thus, although this opinion [9-12] can produce some successful schemes, it is not a general explanation for All-Speed schemes.

Compared with the opinion in [9-12] above, the opinion reflected in Eqs. (50) –

(51) infers that $\Delta U$ itself is reasonable but its coefficient may be incorrect. It then shows great advantages; it can explain all of the uncertainties above and the reason for the success or failing of each scheme in the paper. Therefore, the novelty of this new opinion is that it provides an essential and general understanding of the mechanism of the accuracy problem and then leads to more flexible methods for improving schemes.

The second version is inspired by the second rule, with the pseudo-sound speed terms treated differently in the numerators and the denominators. By modifying the original preconditioned Roe scheme (32) and (37), the new version is obtained as:

$$\delta U = \left[ c' - \frac{1-\theta'}{2} U \frac{U'}{\tilde{c}} - \theta' |U| \right] \frac{\Delta p}{\rho \theta c^2} + \frac{U'}{\tilde{c}} \Delta U, \tag{52}$$

$$\delta p = \frac{U'}{\tilde{c}} \Delta p + \left( c' - |U| + \frac{1-\theta'}{2} U \frac{U'}{\tilde{c}} \right) \rho \Delta U, \tag{53}$$

where $\theta' = \min\left[ M^2, 1 \right]$, $U' = \frac{1}{2}(1+\theta')U$, and $c' = \frac{1}{2}\sqrt{4c^2\theta' + (1-\theta')^2 U^2}$, which are equal to the original in form but are decided by the local parameters instead. Although the global cut-off is still applied to the terms in the denominators for computational stability, it does not affect the accuracy because larger denominators lead to smaller numerical dissipation. Thus, this new preconditioned Roe scheme has better accuracy than the original and will be verified in the following chapter.

## 6. Numerical Experiments

To validate the analysis in Chapter 4 and the newly proposed schemes in Chapter 5, the four numerical cases – lid-driven cavity laminar flow, Euler flow past a cylinder, Euler flow past a turbine blade and a Sod shock tube – are tested in this chapter with the nine schemes discussed above – the Roe, the preconditioned Roe with $K=1$ in Eq. (17) (P-Roe), the low Mach number fix for the Roe (LM-Roe), the Thornber's modified Roe (T-Roe), the All-Speed Roe with the central term (4) (A-Roe-c), the All-Speed Roe with the central term (8) and $c_2 = 0.04$ (A-Roe-p), the All-Speed Roe with time-marching MIM (9) (A-Roe-m), the new All-Speed Roe (50) and (51) (A-Roe-new1), and the new All-Speed Roe (52) and (53) (A-Roe-new2).

The first-order accuracy is adopted for discussion of the schemes themselves, with a second-order result for P-Roe with MUSCL reconstruction (P-Roe-MUSCL) for comparison. The explicit time-marching algorithm is used for computation and the CFL number is taken as 0.8 for the first-order accuracy. The pressure shown in the following figures is defined as an non-dimensional pressure $\bar{p} = \frac{p - p_{min}}{p_{max} - p_{min}}$.

## 6.1 Lid-Driven Cavity Flow

The two-dimensional lid-driven cavity flow problem is a typical low-Mach number test case. The computation is performed with the 160*160 grid points at a Mach number of 0.005 and a Reynolds number of 400.

Fig. 1(a) shows that if the numerical dissipation of the Roe scheme is too large, it prevents it from obtaining a reasonable result, as also shown in Fig. 2. The obvious checkerboard decoupling of A-Roe-c can be observed in Fig. 1(b), though it is much less severe than for the Euler flows in sections 6.2 and 6.3, and the computation remains stable because the physical dissipation damps away the checkerboard to a certain degree, even if does not avoid it fully. Results similar to those in Fig. 1(c) with very weak decoupling can be obtained by the schemes that use the $\mathrm{O}(c^{-1})$ order coefficient for the pressure-derivative term, such as LM-Roe, T-Roe, and A-Roe-new1. As shown in Fig. 1(d), all checkerboard modes can be suppressed by the schemes that use the $\mathrm{O}(c^0)$ order coefficient for the pressure-derivative term, such as P-Roe, A-Roe-p, A-Roe-m, and A-Roe-new2.

Fig. 2 shows the accuracy of the schemes. The solution by Ghia et al. [16] is given as the benchmark, and the second-order solution agrees very well with it. A-Roe-m scheme has slightly better accuracy than A-Roe-p, because it introduces the higher-order pressure-derivative dissipation terms discussed in sections 2.2 and 4.4. The A-Roe-new2 scheme has much better accuracy than P-Roe, because it adopts the local Mach number in the numerators of the coefficients. The result of A-Roe-new2 is very close to that of A-Roe-p. In fact, it is difficult to distinguish among the solutions of A-Roe-new2, A-Roe-p, LM-Roe, T-Roe, and A-Roe-new1,and the latter three are not presented in Fig. 2 for simplicity. Similarly, it is difficult to distinguish between the solutions of A-Roe-m and A-Roe-c, and thus the latter is ignored. This fact also shows that the higher-order pressure-derivative dissipation terms in A-Roe-m produce negligible numerical dissipation from the viewpoint of accuracy, although it is important in terms of the checkerboard.

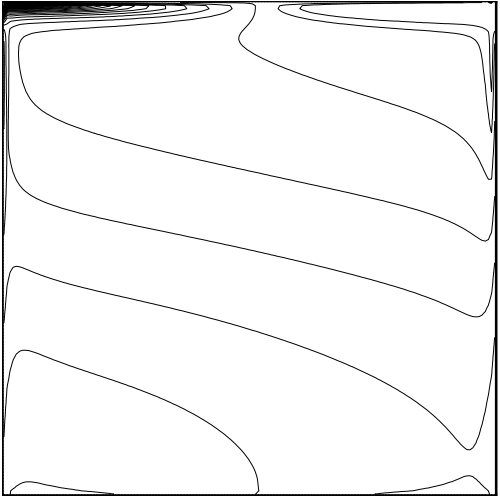 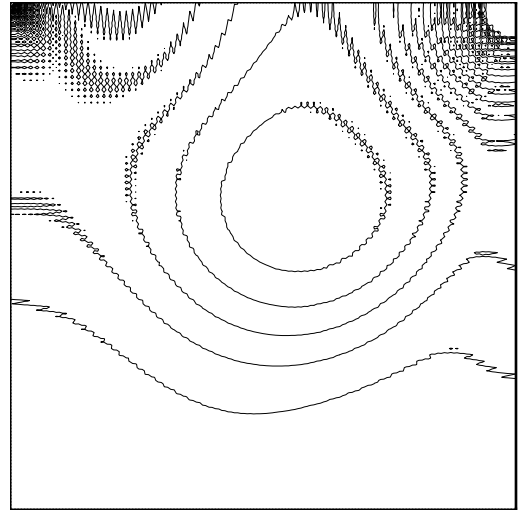

(a) Non-physical solution  (b) Solution with obvious decoupling

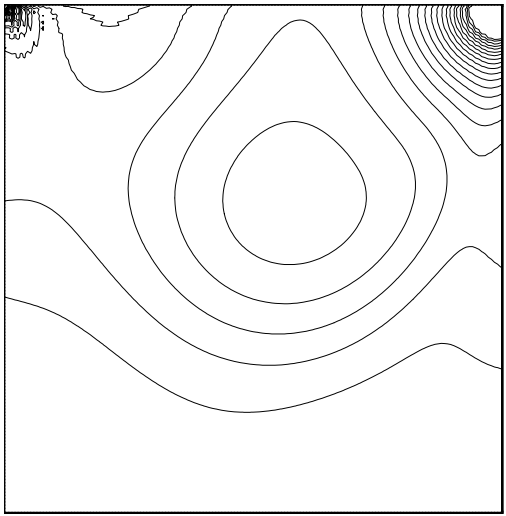 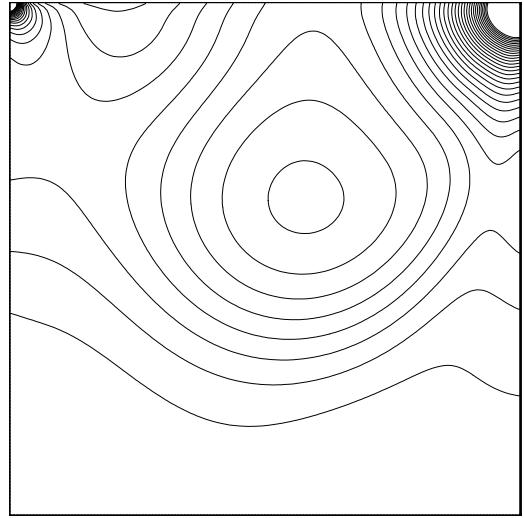

(c) Solution with very weak decoupling  (d) Physical solution

Fig. 1 Pressure contours of lid-driven cavity flow at $M_\infty = 0.005$ and $Re = 400$

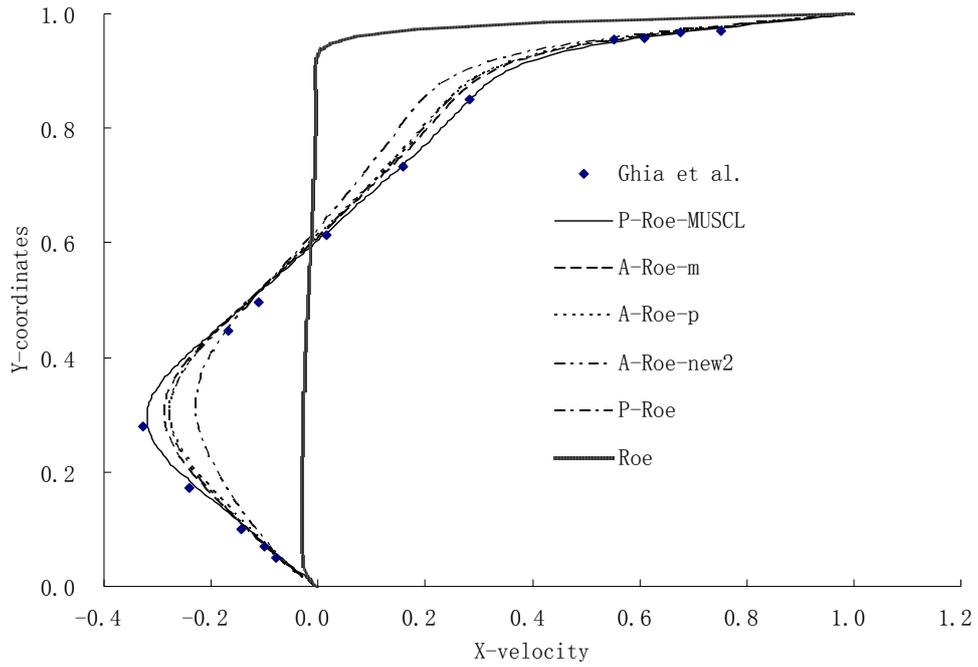

Fig. 2 Computed X-direction velocity distribution along the Y-direction geometric

centreline

## 6.2 Inviscid Flow Around a Cylinder

The two-dimensional Euler flow past a cylinder is another typical low-Mach number test case. The computation is performed with an inflow Mach number of 0.01 and the 72*100 O-type grid points along the circumference and radius, respectively.

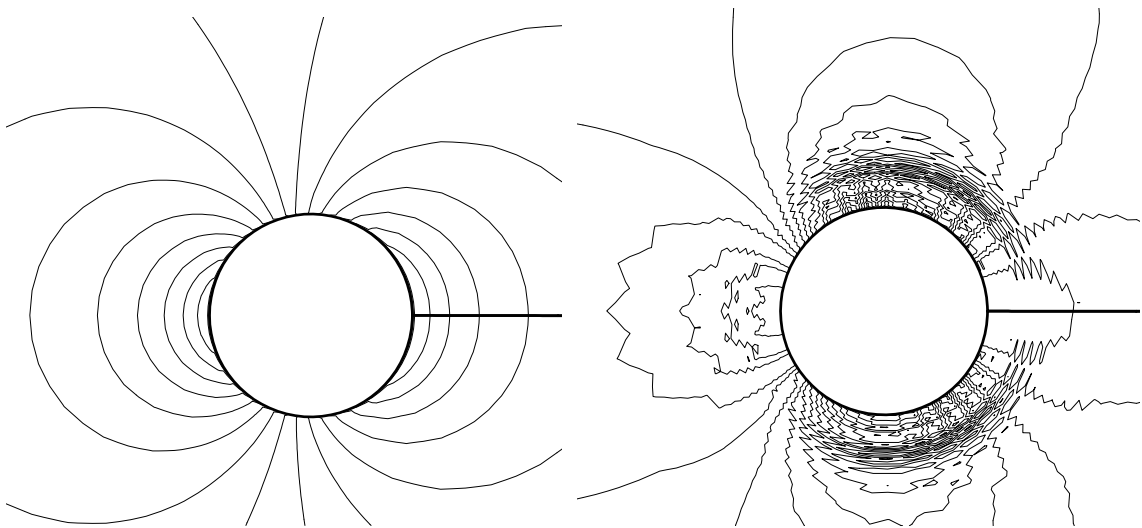

(a) Solution similar to the Stokes flow    (b) Solution with full decoupling

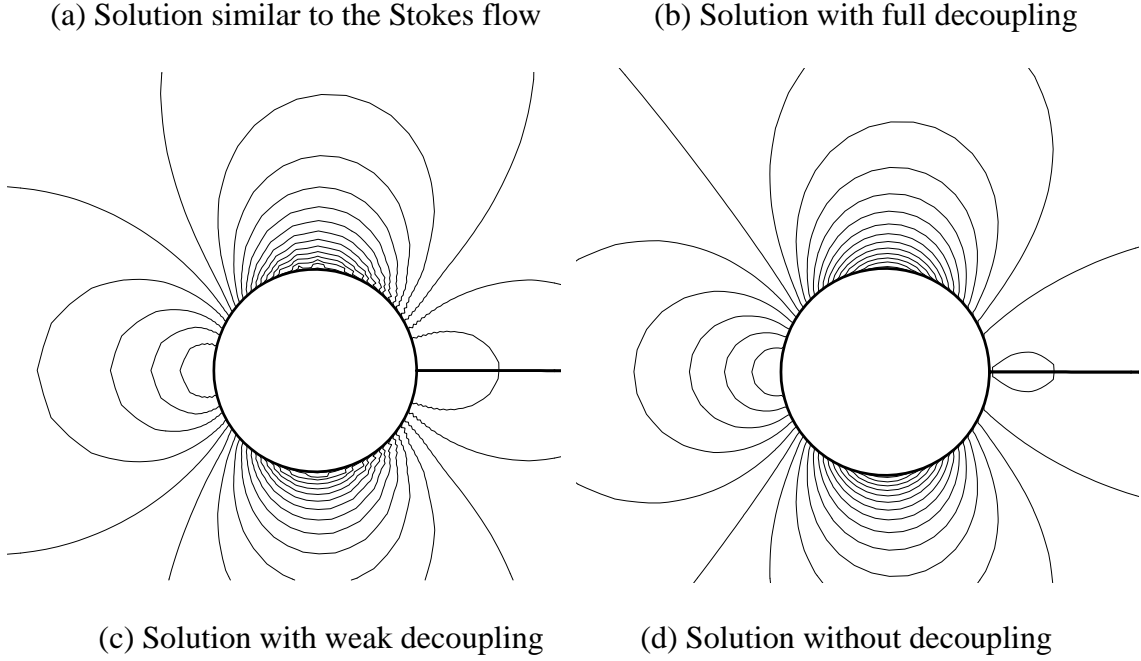

(c) Solution with weak decoupling    (d) Solution without decoupling

Fig. 3 Pressure contours of the inviscid flow around a cylinder at $M_\infty = 0.01$

The Roe scheme produces the result shown in Fig. 1(a), resembling the Stokes flow, which is explained in Ref. [10]. With A-Roe-c, a severe pressure checkerboard quickly appears, as shown in Fig. 3(b), even if the good solution in Fig. 3(d) is used as the initial field. The effect is much more severe than for the viscous simulation in Fig. 1(b) and increases over time, leading to instability. The schemes with the $\mathrm{O}(c^{-1})$ order coefficient for the pressure-derivative term, such as LM-Roe, T-Roe and A-Roe-new1, produce weak checkerboard decoupling near the cylinder as shown in Fig. 3(c). However, the schemes with the $\mathrm{O}(c^0)$ order coefficient for the pressure-derivative term, such as A-Roe-p, A-Roe-m, A-Roe-new2, and P-Roe, can fully suppress the checkerboard as shown in Fig. 3(d).

The following figures provide a quantitative comparison to further verify the accuracy of the schemes. Besides using P-Roe and P-Roe-MUSCL as the benchmark, only A-Roe-new1 and A-Roe-new2 are compared, because the accuracy of other schemes has been validated in the relevant papers. Fig. 4 gives the non-dimensional pressure distribution on the cylinder surface and shows that A-Roe-new1 and A-Roe-new2 have better accuracy than P-Roe. Following Ref. [1], Fig. 5 displays the pressure fluctuations Ind($p$) = ($P_{max}$-$P_{min}$) / $P_{max}$ versus the inlet Mach number, which perfectly agrees with the theoretical asymptotic predictions: the pressure fluctuations scale exactly with $M_*^2$.

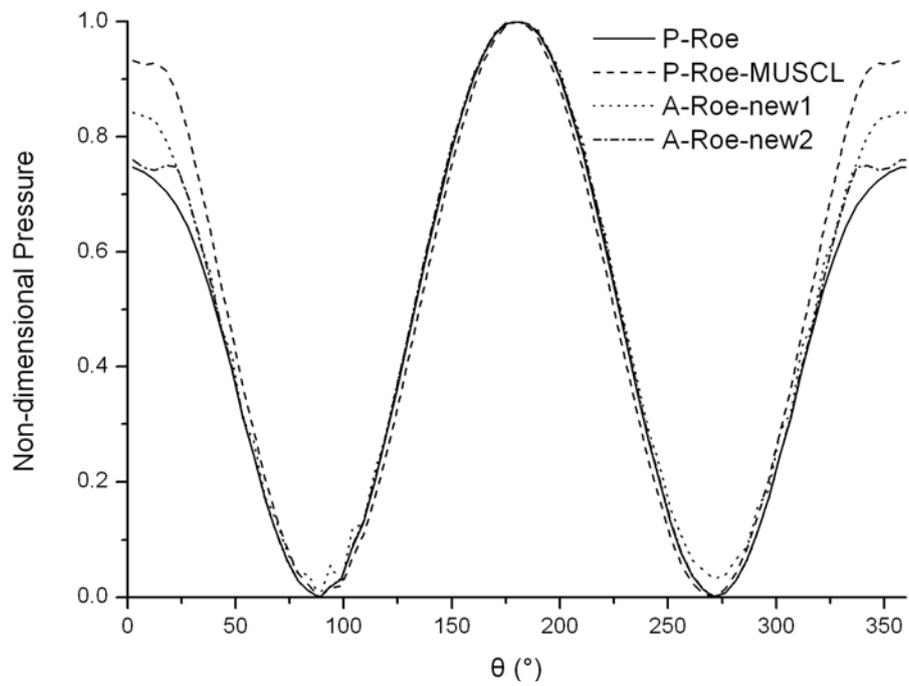

Fig. 4 Pressure on the cylinder surface at $M_\infty = 0.01$

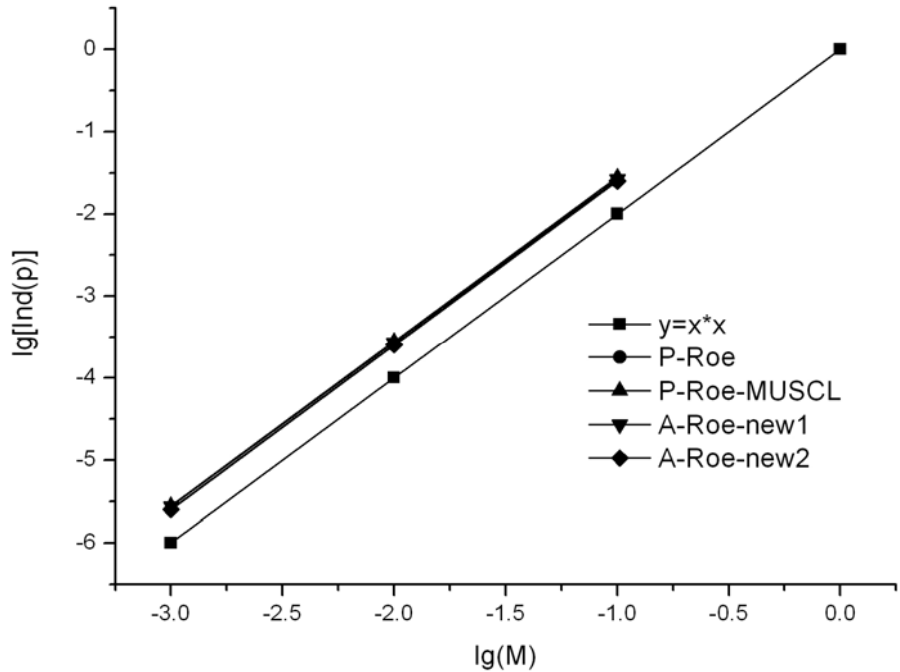

Fig. 5 Pressure fluctuations vs inflow Mach number

## 6.3 Inviscid Flow Around a Turbine Blade

The Euler flow past a high-load turbine blade (T106) row is simulated with the inlet Mach number 0.001, with the 40*98 H-type grid points in the azimuthal and streamline directions.

Similar to the inviscid flow around the cylinder, the Roe scheme obtains a nonphysical solution as shown in Fig. 6(a). A-Roe-c suffers from the severe checkerboard problem, as shown in Fig. 6(b), and is unstable. Schemes LM-Roe, T-Roe, and A-Roe-new1 suffer form the obvious checkerboard problem, as shown in Fig. 6(c), but remain stable, and schemes A-Roe-p, A-Roe-m, A-Roe-new2, and P-Roe can suppress the checkerboard fully as shown in Fig. 6(d). The quantitative comparisons in Fig. 7 and Fig. 8 also show that A-Roe-new1 and A-Roe-new2 work well.

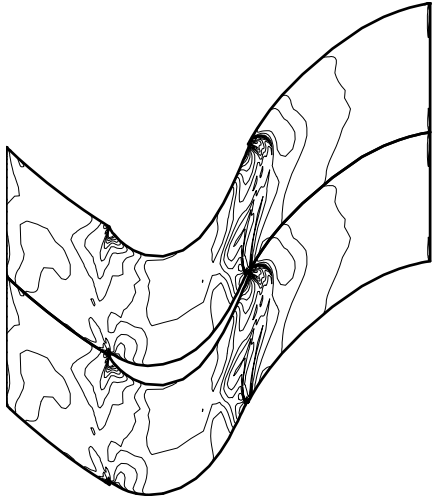
(a) Non-physical solution

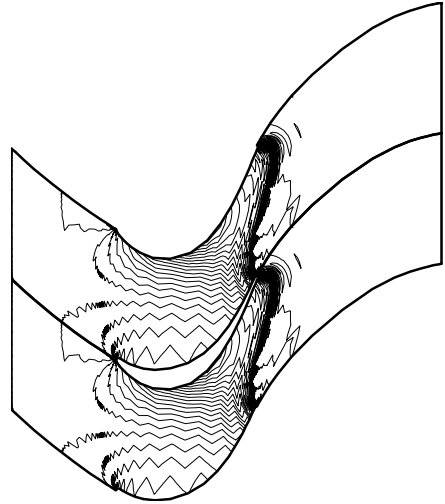
(b) Solution with full decoupling

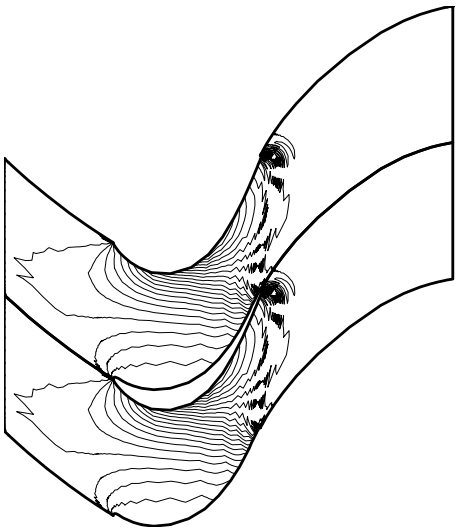
(c) Solution with weak decoupling

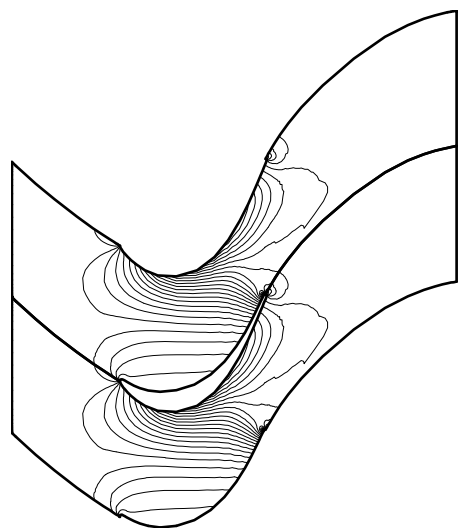
(d) Solution without decoupling

Fig. 6 Pressure contours of inviscid flow around the turbine blade row at $M_{in} = 0.001$

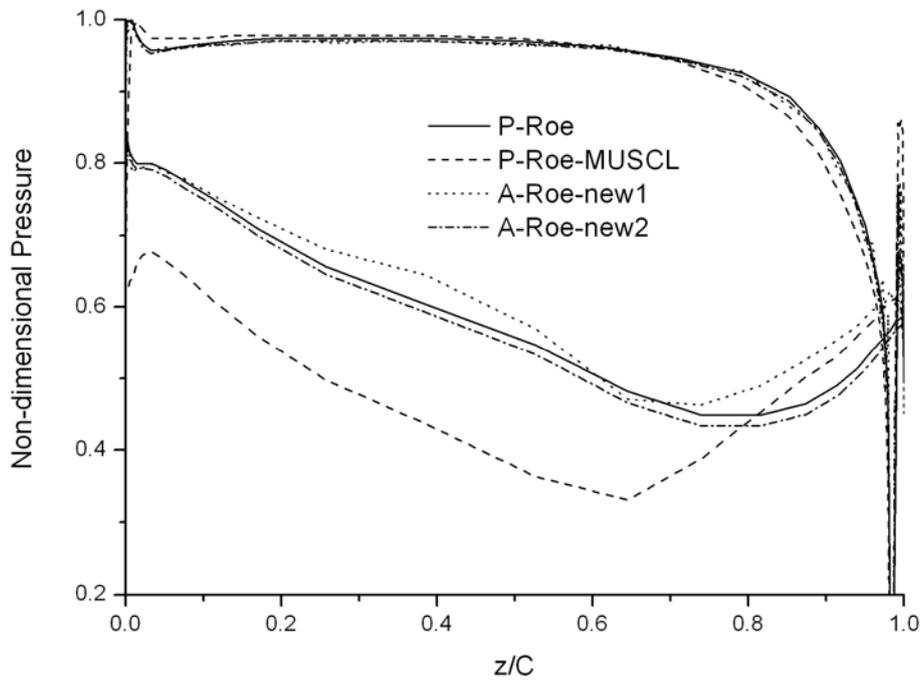

Fig. 7 Pressure on the blade surface at $M_{in} = 0.001$

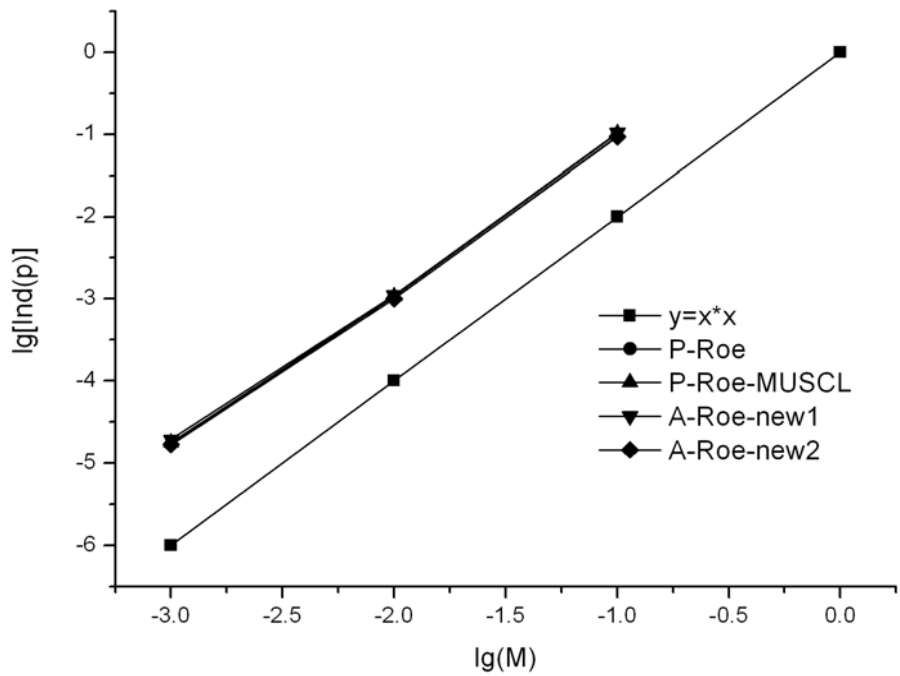

Fig. 8 Pressure fluctuations vs inflow Mach number

### 6.4 Sod Shock Tube

To demonstrate that the proposed new Roe-type schemes can handle shocks, the

Sod shock tube test case is adopted with the 200 cells and the initial conditions as follows:

$$(\rho, u, p)_L = (1, 0, 1) \text{ if } x < 0.5 \text{ and } (\rho, u, p)_R = (0.125, 0, 0.1) \text{ if } x > 0.5.$$

The results shown in Fig. 9 and Fig. 10 are taken at $t = 0.2$. Compared with the original Roe scheme, the new Roe-type schemes produce the slightly steeper gradients in the profiles, but A-Roe-new1 has a solution with little jump near $x = 0.5$.

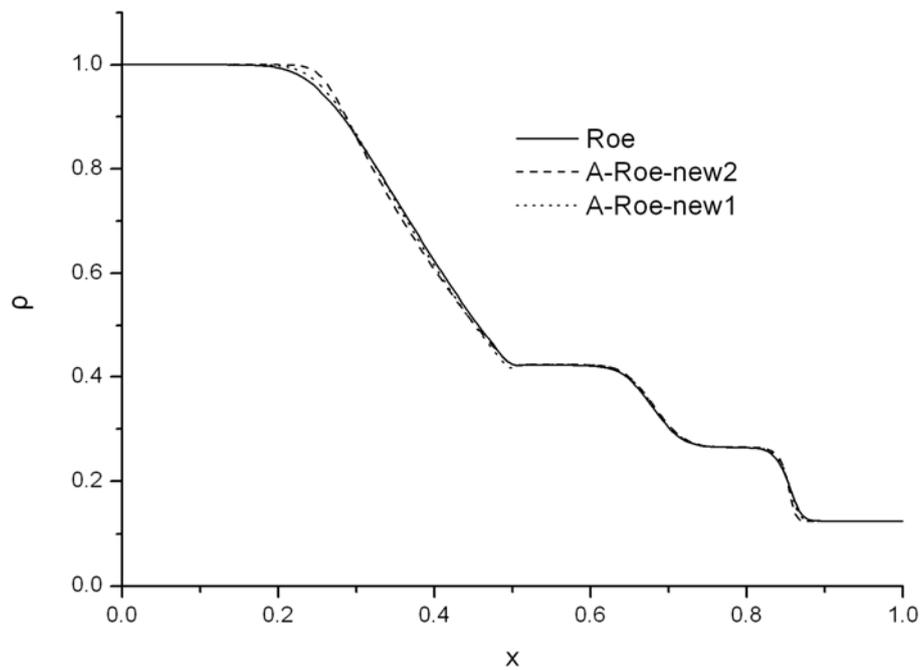

Fig. 9 Density change of the Sod shock tube test

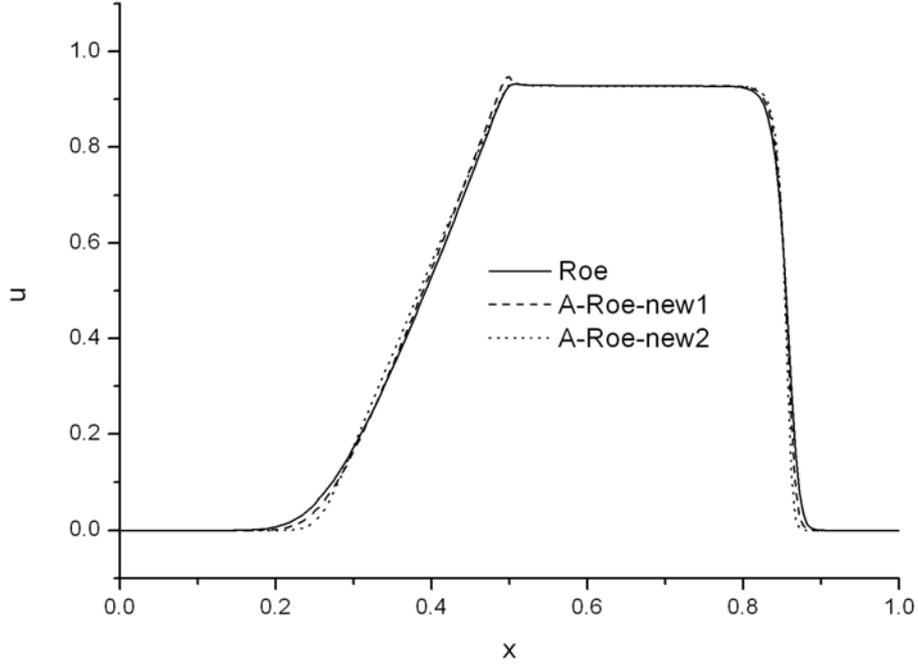

Fig. 10 Velocity change of the Sod shock tube test

## 7. Conclusions

The Roe scheme and its four modified schemes for all-speed flows are compared and analyzed. An explanation is provided for why these Roe-type schemes succeed or fail with the accuracy and checkerboard problems. It can be concluded that the accuracy of the schemes is controlled by the order of $c$ in the coefficient of the term $\Delta U$ in $\delta p$ in the momentum equation, and that the checkerboard is controlled by the order of $c$ in the coefficient of the term $\Delta p$ in $\delta U$, especially in the continuity equation. Based on this understanding, the following rules are proposed for constructing more satisfactory schemes for all-speed flows.

(1) The coefficient of the velocity-derivative dissipation term should be of the order $O(c^0)$ or lower for accuracy. The coefficient of the pressure-derivative dissipation term should also be of the order $O(c^0)$ to suppress the checkerboard completely with the divergence constraint of $\boldsymbol{u}^0$ unsatisfied, or of the order $\mathrm{O}(c^{-1})$ with a weak checkerboard, or of in-between.

(2) For accuracy, it is enough to multiply the sound speed term in the numerator with the function $f(M)$, which only relates to the local variables. If the sound speed in the denominator changes, a global cut-off strategy is necessary, which does not really decrease accuracy because the numerical dissipation reduces as the value of the denominator increases. For example, the preconditioned Roe scheme suffers from the global cut-off because it treats the sound speed terms in both the numerator and the denominator equally. A more satisfactory All-Speed scheme should treat these terms differently.

Two novel schemes are proposed as examples based on these two rules. Numerical experiments show that the behaviour of all the schemes discussed herein accord with the theoretical prediction.

## Acknowledgments

This work is supported by Project 50806037 of the National Natural Science Foundation of China and Project 2007CB210105 of the 973 Program. We would like to thank Prof. Wang Cuncheng for polishing the paper.